# A course on Tug of War games with random noise

Introduction and basic constructions

Marta Lewicka

# Preface

The goal of these Course Notes is to present a systematic overview of the basic constructions and results pertaining to the recently emerged field of Tug of War games, as seen from an analyst's perspective. To a large extent, this book represents the author's own study itinerary, aiming at precision and completeness of a classroom text in an upper undergraduate to graduate level course.

This book was originally planned as a joint project between Marta Lewicka (University of Pittsburgh) and Yuval Peres (then Microsoft Research). Due to an unforeseen turn of events, neither the collaboration nor the execution of the project in the priorly conceived forms, could have been pursued.

The author wishes to dedicate this book to all women in mathematics, with admiration and encouragement. The publishing profit will be donated to the Association for Women in Mathematics.

<div align="right">

Marta Lewicka,
Pittsburgh, October 2019.

</div>

# 2

# The linear case: random walk and harmonic functions

In this Chapter we present the basic relation between the *linear potential theory* and *random walks*. This fundamental connection, developed by Ito, Doob, Lévy and others, relies on the observation that harmonic functions and martingales share a common cancellation property, expressed via *mean value properties*. It turns out that, with appropriate modifications, a similar observation and approach can be applied also in the nonlinear case, which is of main interest in these Course Notes. Thus, the present Chapter serves as a stepping stone towards gaining familiarity with more complex constructions of Chapters 3-6.

After recalling the equivalent defining properties of *harmonic functions* in Section 2.1, in Section 2.2 we introduce the *ball walk*. This is an infinite horizon discrete process, in which at each step the particle, initially placed at some point $x_0$ in the open, bounded domain $\mathcal{D} \subset \mathbb{R}^N$, is randomly advanced to a new position, uniformly distributed within the following open ball: centered at the current placement, and with radius equal to the minimum of the parameter $\epsilon$ and the distance from the boundary $\partial \mathcal{D}$. With probability one, such process accumulates on $\partial \mathcal{D}$ and $u^\epsilon(x_0)$ is then defined as the expected value of the given boundary data $F$ at the process limiting position. Each function $u^\epsilon$ is harmonic, and we show in Sections 2.3 and 2.4, that if $\partial \mathcal{D}$ is *regular*, then each $u^\epsilon$ coincides with the unique *harmonic extension* of $F$ in $\mathcal{D}$. One sufficient condition for regularity is the *exterior cone condition*, as proved in Section 2.5.

Our discussion and proofs are elementary, requiring only a basic knowledge of probabilistic concepts, such as: probability spaces, martingales and Doob's theorem. For convenience of the reader, these are gathered in Appendix A. The slightly more advanced material which may be skipped at first reading, is based on the Potential Theoretic and the Brownian motion arguments from, respectively, Appendix C and Appendix B. Both approaches allow to deduce that functions in the family $\{u^\epsilon\}_{\epsilon \in (0,1)}$ are one and the same function, regardless of the regularity of $\partial \mathcal{D}$. This fact is obtained first in Section 2.6* by proving that $u^\epsilon$ coincide with the *Perron solution* of the Dirichlet problem for boundary data $F$. The same follows in Section 2.7* by checking that the ball walk consists of discrete realisations along the Brownian motion trajectories, to the effect that $u^\epsilon$ equal the *Brownian motion harmonic extension* of $F$.

Thus, the three classical approaches to finding the harmonic extension by:

(i) evaluating the expectation of the values of the (discrete) ball walk at its limiting infinite horison boundary position;
(ii) taking infima/suprema of super- and sub-harmonic functions obeying comparison with the boundary data;
(iii) evaluating the expectation of the values of the (continuous) Brownian motion at exiting the domain;

are shown to naturally coincide when $F$ is continuous.

## 2.1 The Laplace equation and harmonic functions

Among the most important of all PDEs is the *Laplace equation*. In this Section we briefly recall the relevant definitions and notation; for the proofs and a review of basic properties we refer to Section C.3 in Appendix C.

Let $\mathcal{D} \subset \mathbb{R}^N$ be an open, bounded, connected set. The Euler-Lagrange equation for critical points of the following quadratic energy functional:

$$\mathcal{I}_2(u) = \int_{\mathcal{D}} |\nabla u(x)|^2 \, \mathrm{d}x$$

is expressed by the second order partial differential equation:

$$\Delta u \doteq \sum_{i=1}^{N} \frac{\partial^2 u}{(\partial x_i)^2} = 0 \quad \text{in } \mathcal{D},$$

whose solutions are called *harmonic functions*. The operator $\Delta$ is defined in the classical sense only for $C^2$ functions $u$, however a remarkable property of harmonicity is that it can be equivalently characterised via *mean value properties* that do not require $u$ to be even continuous. At the same time, harmonic functions are automatically smooth. More precisely, the following conditions are equivalent (the proof will be recalled in Section C.3):

(i) A locally bounded, Borel function $u : \mathcal{D} \to \mathbb{R}$ satisfies the *mean value property on balls*:

$$u(x) = \fint_{B_r(x)} u(y) \, \mathrm{d}y \qquad \text{for all } \bar{B}_r(x) \subset \mathcal{D}.$$

(ii) A locally bounded, Borel function $u : \mathcal{D} \to \mathbb{R}$ satisfies for each $x \in \mathcal{D}$ and almost every $r \in (0, \text{dist}(x, \partial\mathcal{D}))$ the *mean value property on spheres*:

$$u(x) = \fint_{\partial B_r(x)} u(y) \, d\sigma^{N-1}(y).$$

(iii) The function $u$ is smooth: $u \in C^\infty(\mathcal{D})$, and there holds $\Delta u = 0$ in $\mathcal{D}$.

We also remark at this point that, Taylor expanding any function $u \in C^2(\mathcal{D})$ and averaging term by term on $\bar{B}_\epsilon(x) \subset \mathcal{D}$, leads to the *mean value expansion*, also called in what follows the *averaging principle*:

$$\fint_{B_\epsilon(x)} u(y) \, dy = u(x) + \frac{\epsilon^2}{2(N+2)} \Delta u(x) + o(\epsilon^2) \qquad \text{as } \epsilon \to 0+, \qquad (2.1)$$

which in fact is consistent with interpreting $\Delta u$ as the (second order) error from harmonicity. This point of view is central to developing the probabilistic interpretation of the general $p$-Laplace equations, which is the goal of these Course Notes. While we will not need (2.1) in order to construct the random walk and derive its connection to the Laplace equation $\Delta$ in the linear setting $p = 2$ studied in this Chapter, it is beneficial to keep in mind that the mean value property in (i) may be actually "guessed" from the expansion (2.1).

Throughout next Chapters, more general averaging principles will be proved (in Sections 3.2, 4.1 and 6.1), informing the mean value properties that characterise, in the asymptotic sense, zeroes of the nonlinear operators $\Delta_p$ at any $p \in (1, \infty)$, and ultimately leading to the Tug of War games with random noise.

## 2.2 The ball walk

In this Section we construct the discrete stochastic process whose value will be shown to equal the harmonic function with prescribed boundary values.

The probability space of the ball walk process is defined as follows. Consider $(\Omega_1, \mathcal{F}_1, \mathbb{P}_1)$, where $\Omega_1$ is the unit ball $B_1(0) \subset \mathbb{R}^N$, the $\sigma$-algebra $\mathcal{F}_1$ consists of Borel subsets of $\Omega_1$, and $\mathbb{P}_1$ is the normalised Lebesgue measure:

$$\mathbb{P}_1(D) = \frac{|D|}{|B_1(0)|} \qquad \text{for all } D \in \mathcal{F}_1,$$

For any $n \in \mathbb{N}$, we denote by $\Omega_n = (\Omega_1)^n$ the Cartesian product of $n$ copies of $\Omega_1$, and by $(\Omega_n, \mathcal{F}_n, \mathbb{P}_n)$ the corresponding product probability space. Further, the countable product $(\Omega, \mathcal{F}, \mathbb{P})$ is defined as in Theorem A.12 on:

$$\Omega \doteq (\Omega_1)^\mathbb{N} = \prod_{i=1}^\infty \Omega_1 = \Big\{ \omega = \{w_i\}_{i=1}^\infty; \ w_i \in B_1(0) \ \text{ for all } i \in \mathbb{N} \Big\}.$$

We identify each $\sigma$-algebra $\mathcal{F}_n$ with the sub-$\sigma$-algebra of $\mathcal{F}$ consisting of sets of the form $F \times \prod_{i=n+1}^{\infty} \Omega_1$ for all $F \in \mathcal{F}_n$. Note that $\{\mathcal{F}_n\}_{n=0}^{\infty}$ where $\mathcal{F}_0 = \{\emptyset, \Omega\}$, is a filtration of $\mathcal{F}$ and that $\mathcal{F}$ is the smallest $\sigma$-algebra containing $\bigcup_{n=0}^{\infty} \mathcal{F}_n$.

**Definition 2.1** Let $\mathcal{D} \subset \mathbb{R}^N$ be an open, bounded, connected set. The *ball walk* for $\epsilon \in (0, 1)$ and $x_0 \in \mathcal{D}$, is recursively defined (see Figure 2.1) as the following sequence of random variables $\{X_n^{\epsilon, x_0} : \Omega \to \mathcal{D}\}_{n=0}^{\infty}$:

$$X_0^{\epsilon, x_0} \equiv x_0,$$
$$X_n^{\epsilon, x_0}(w_1, \ldots, w_n) = X_{n-1}^{\epsilon, x_0}(w_1, \ldots, w_{n-1}) + (\epsilon \wedge \operatorname{dist}(X_{n-1}^{\epsilon, x_0}, \partial \mathcal{D})) w_n \quad (2.2)$$
$$\text{for all } n \geq 1 \text{ and all } (w_1, \ldots w_n) \in \Omega_n.$$

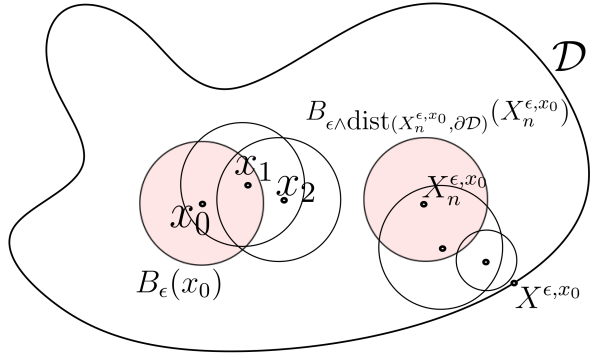

Figure 2.1 The ball walk and the process $\{X_n^{\epsilon, x_0}\}_{n=0}^{\infty}$ in (2.2).

We will often write: $x_n = X_n^{\epsilon, x_0}(w_1, \ldots, w_n)$. Intuitively, $\{x_n\}_{n=0}^{\infty}$ describe the consecutive positions of a particle initially placed at $x_0 \in \mathcal{D}$, along a discrete path consisting of a succession of random steps of magnitude at most $\epsilon$. The size of steps decreases as the particle approaches the boundary $\partial \mathcal{D}$. The position $x_n \in \mathcal{D}$ is obtained from $x_{n-1}$ by sampling uniformly on the open ball $B_{\epsilon \wedge \operatorname{dist}(x_{n-1}, \partial \mathcal{D})}(x_{n-1})$. It is clear that each random variable $X_n^{\epsilon, x_0} : \Omega \to \mathbb{R}^N$ is $\mathcal{F}_n$-measurable and that it depends only on the previous position $x_{n-1}$, its distance from $\partial \mathcal{D}$ and the current random outcome $w_n \in \Omega_1$.

**Lemma 2.2** *In the above context, the sequence $\{X_n^{\epsilon, x_0}\}_{n=0}^{\infty}$ is a martingale relative to the filtration $\{\mathcal{F}_n\}_{n=0}^{\infty}$, namely:*

$$\mathbb{E}(X_n^{\epsilon, x_0} \mid \mathcal{F}_{n-1}) = X_{n-1}^{\epsilon, x_0} \quad \mathbb{P} - a.s \quad \text{for all } n \geq 1.$$

*Moreover, there exists a random variable $X^{\epsilon,x_0} : \Omega \to \partial\mathcal{D}$ such that:*

$$\lim_{n\to\infty} X_n^{\epsilon,x_0} = X^{\epsilon,x_0} \qquad \mathbb{P}-a.s. \tag{2.3}$$

*Proof* **1.** Since the sequence $\{X_n^{\epsilon,x_0}\}_{n=0}^{\infty}$ is bounded in view of boundedness of $\mathcal{D}$, Theorem A.38 will yield convergence in (2.3) provided we check the martingale property. Indeed it follows that (see Lemma A.17):

$$\mathbb{E}(X_n^{\epsilon,x_0} \mid \mathcal{F}_{n-1})(w_1,\ldots,w_{n-1}) = \int_{\Omega_1} X_n^{\epsilon,x_0}(w_1,\ldots,w_n)\,d\mathbb{P}_1(w_n)$$

$$= x_{n-1} + (\epsilon \wedge \text{dist}(x_{n-1}, \partial\mathcal{D})) \int_{\Omega_1} w_n\,d\mathbb{P}_1(w_n)$$

$$= X_{n-1}^{\epsilon,x_0}(w_1,\ldots,w_{n-1}) \quad \text{for } \mathbb{P}_{n-1}\text{-a.e. } (w_1,\ldots,w_{n-1}) \in \Omega_{n-1}.$$

**2.** It remains to prove that the limiting random variable $X^{\epsilon,x_0} : \Omega \to \bar{\mathcal{D}}$ satisfies $\mathbb{P}$-a.s. the boundary accumulation property: $X^{\epsilon,x_0} \in \partial\mathcal{D}$. Observe that:

$$\{\lim_{n\to\infty} X_n^{\epsilon,x_0} = X^{\epsilon,x_0}\} \cap \{X^{\epsilon,x_0} \in \mathcal{D}\} \subset \bigcup_{n\in\mathbb{N},\ \delta\in(0,\epsilon)\cap\mathbb{Q}} A(n,\delta), \tag{2.4}$$

where $A(n,\delta) = \{\text{dist}(X_i^{\epsilon,x_0}, \partial\mathcal{D}) \geq \delta \text{ and } |X_{i+1}^{\epsilon,x_0} - X_i^{\epsilon,x_0}| \leq \frac{\delta}{2} \text{ for all } i \geq n\}$. Then:

$$A(n,\delta) \subset \{\omega \in \Omega;\ |w_i| \leq \frac{1}{2} \text{ for all } i > n\}.$$

Indeed, if $\omega = \{w_i\}_{i=1}^{\infty} \in A(n,\delta)$ with $\delta < \epsilon$, it follows that:

$$\frac{\delta}{2} \geq |X_{i+1}^{\epsilon,x_0}(\omega) - X_i^{\epsilon,x_0}(\omega)| = (\epsilon \wedge \text{dist}(X_i^{\epsilon,x_0}(\omega), \partial\mathcal{D}))\,|w_{i+1}| \geq (\epsilon \wedge \delta)|w_{i+1}| = \delta|w_{i+1}|,$$

which implies that $|w_{i+1}| \leq \frac{1}{2}$ for all $i \geq n$. Concluding:

$$\mathbb{P}(A(n,\delta)) \leq \lim_{i\to\infty} \mathbb{P}_1(B_{\frac{1}{2}}(0))^{i-n} = 0 \quad \text{for all } n \in \mathbb{N} \text{ and all } \delta \in (0,\epsilon).$$

Hence, the event in the left hand side of (2.4) has probability 0. $\square$

Given now a continuous function $F : \partial\mathcal{D} \to \mathbb{R}$, define:

$$u^{\epsilon}(x_0) \doteq \mathbb{E}[F \circ X^{\epsilon,x_0}] = \int_{\Omega} F \circ X^{\epsilon,x_0}\,d\mathbb{P}. \tag{2.5}$$

Note that the above construction obeys the *comparison principle*. Namely, if $F, \bar{F} : \partial\mathcal{D} \to \mathbb{R}$ are two continuous functions such that $F \leq \bar{F}$ on $\partial\mathcal{D}$, then the corresponding $u^{\epsilon}$ and $\bar{u}^{\epsilon}$ satisfy: $u^{\epsilon} \leq \bar{u}^{\epsilon}$ in $\mathcal{D}$.

**Remark 2.3** It is useful to view the boundary function $F$ as the restriction on $\partial \mathcal{D}$ of some continuous $F : \bar{\mathcal{D}} \to \mathbb{R}$, see Exercise 2.7 (i). Then we may write:

$$u^\epsilon(x_0) = \lim_{n \to \infty} \int_\Omega F \circ X_n^{\epsilon, x_0} \, d\mathbb{P}. \tag{2.6}$$

Since for each $n \geq 0$ the function $F \circ X_n^{\epsilon, x_0}$ is jointly Borel-regular in the variables $x_0 \in \mathcal{D}$ and $\omega \in \Omega_n$, it follows by Theorem A.11 that $x_0 \mapsto \mathbb{E}[F \circ X_n^{\epsilon, x_0}]$ is Borel-regular. Consequently, $u^\epsilon : \mathcal{D} \to \mathbb{R}$ is also Borel.

In what follows, we will denote the average $\mathcal{A}_\delta u$ of an integrable function $u : \mathcal{D} \to \mathbb{R}$ on a ball $B_\delta(x) \subset \mathcal{D}$ by:

$$\mathcal{A}_\delta u(x) \doteq \fint_{B_\delta(x)} u(y) \, dy.$$

Directly from Definition 2.1 and (2.5) we conclude the satisfaction of the mean value property for each $u^\epsilon$ on the sampling balls from (2.2):

---

**Theorem 2.4** *Let $\mathcal{D} \subset \mathbb{R}^N$ be open, bounded, connected, and let $F : \partial \mathcal{D} \to \mathbb{R}$ be continuous. Then, the function $u^\epsilon : \mathcal{D} \to \mathbb{R}$ defined in (2.5) and equivalently in (2.6), is continuous and satisfies:*

$$u^\epsilon(x) = \mathcal{A}_{\epsilon \wedge dist(x, \partial \mathcal{D})} u^\epsilon(x) \quad \text{for all } x \in \mathcal{D}.$$

---

*Proof* Fix $\epsilon \in (0, 1)$ and $x_0 \in \mathcal{D}$. For each $n \geq 2$ we view $(\Omega_n, \mathcal{F}_n, \mathbb{P}_n)$ as the product of probability spaces $(\Omega_1, \mathcal{F}_1, \mathbb{P}_1)$ and $(\Omega_{n-1}, \mathcal{F}_{n-1}, \mathbb{P}_{n-1})$. Applying Fubini's Theorem (Theorem A.11), we get:

$$\mathbb{E}[F \circ X_n^{\epsilon, x_0}] = \int_{\Omega_1} \int_{\Omega_{n-1}} (F \circ X_n^{\epsilon, x_0})(w_1, \ldots, w_n) \, d\mathbb{P}_{n-1}(w_2, \ldots, w_n) \, d\mathbb{P}_1(w_1)$$

$$= \int_{\Omega_1} \mathbb{E}[F \circ X_{n-1}^{\epsilon, X_1^{\epsilon, x_0}(w_1)}] \, d\mathbb{P}_1(w_1),$$

where $F : \bar{\mathcal{D}} \to \mathbb{R}$ is some continuous extension of its given values on $\partial \mathcal{D}$, as in (2.6). Passing to the limit with $n \to \infty$ and changing variables, we obtain:

$$u^\epsilon(x_0) = \int_{\Omega_1} u^\epsilon(X_1^{\epsilon, x_0}(w_1)) \, d\mathbb{P}_1(w_1)$$

$$= \int_{\Omega_1} u^\epsilon(x_0 + (\epsilon \wedge \text{dist}(x_0, \partial \mathcal{D})) w_1) \, d\mathbb{P}_1(w_1)$$

$$= \fint_{B_{\epsilon \wedge \text{dist}(x_0, \mathcal{D})}(x_0)} u^\epsilon(y) \, dy.$$

Continuity of $u^\epsilon$ follows directly from the averaging formula and we leave it as an exercise (see Exercise 2.7 (ii)). □

The next two statements imply uniqueness of classical solutions to the boundary value problem for the Laplacian. The same property, in the basic analytical setting that we review in Section C.3, follows via the maximum principle.

**Corollary 2.5** *Let $\epsilon \in (0,1)$, $x_0 \in \mathcal{D}$ and let $\{X_n^{\epsilon,x_0}\}_{n=0}^\infty$ be as in (2.2). In the setting of Theorem 2.4, the sequence $\{u^\epsilon \circ X_n^{\epsilon,x_0}\}_{n=0}^\infty$ is a martingale relative to the filtration $\{\mathcal{F}_n\}_{n=0}^\infty$.*

*Proof* Indeed, Lemma A.17 yields for all $n \geq 1$:

$$\mathbb{E}(u^\epsilon \circ X_n^{\epsilon,x_0} \mid \mathcal{F}_{n-1})(w_1,\ldots,w_{n-1}) = \int_{\Omega_1} (u^\epsilon \circ X_n^{\epsilon,x_0})(w_1,\ldots,w_n) \, \mathrm{d}\mathbb{P}_1(w_n)$$

$$= \int_{\Omega_1} u^\epsilon\Big(X_{n-1}^{\epsilon,x_0}(w_1,\ldots,w_{n-1}) + (\epsilon \wedge \mathrm{dist}(x_{n-1},\partial\mathcal{D}))\, w_n\Big) \, \mathrm{d}\mathbb{P}_1(w_n)$$

$$= \fint_{B_{\epsilon \wedge \mathrm{dist}(x_{n-1},\mathcal{D})}(x_{n-1})} u^\epsilon(y) \, \mathrm{d}y = (u^\epsilon \circ X_{n-1}^{x_0})(w_1,\ldots,w_{n-1}),$$
(2.7)

valid for $\mathbb{P}_{n-1}$-a.e. $(w_1,\ldots,w_{n-1}) \in \Omega_{n-1}$. □

---

**Lemma 2.6** *In the setting of Theorem 2.4, assume that $u \in C(\bar{\mathcal{D}})$ solves:*

$$\Delta u = 0 \quad \text{in } \mathcal{D}, \qquad u = F \quad \text{on } \partial\mathcal{D}. \tag{2.8}$$

*Then $u^\epsilon = u$ for all $\epsilon \in (0,1)$. In particular, (2.8) has at most one solution.*

---

*Proof* We first claim that given $x_0 \in \mathcal{D}$ and $\epsilon \in (0,1)$, the sequence $\{u \circ X_n^{\epsilon,x_0}\}_{n=0}^\infty$ is a martingale relative to $\{\mathcal{F}_n\}_{n=0}^\infty$. This property follows exactly as in (2.7), where $u^\epsilon$ is now replaced by $u$ and where the mean value property for harmonic functions (C.8) is used instead of the single-radius averaging formula of Theorem 2.4. Consequently, we get:

$$u(x_0) = \mathbb{E}[u \circ X_0^{\epsilon,x_0}] = \mathbb{E}[u \circ X_n^{\epsilon,x_0}] \quad \text{for all } n \geq 0.$$

Since the right hand side above converges to $u^\epsilon(x_0)$ with $n \to \infty$, it follows that $u(x_0) = u^\epsilon(x_0)$. To prove the second claim, recall that $u^\epsilon(x_0)$ depends only on the boundary values $u_{|\partial\mathcal{D}} = F$ and not on their extension $u$ on $\bar{\mathcal{D}}$. This yields uniqueness of the harmonic extension in (2.8). □

We finally remark that the mean value property stated in Theorem 2.4 suffices to conclude that each $u^\epsilon$ is harmonic (see Section C.3). One can also show that all functions in the family $\{u^\epsilon\}_{\epsilon\in(0,1)}$ are the same, even in the absence of the classical harmonic extension $u$ satisfying (2.8). This general result will be given two independent proofs in Sections 2.6* and 2.7*. In the next Section, we provide an elementary proof in domains that are sufficiently regular. An entirely similar strategy, based on showing the uniform convergence of $\{u^\epsilon\}_{\epsilon\to 0}$ in $\bar{\mathcal{D}}$ and analyzing its limit, will be adopted in Chapters 3-6 for the $p$-harmonic case, $p \in (1, \infty)$, in the context of Tug of War with noise.

**Exercise 2.7**  (i) Let $F : A \to \mathbb{R}$ be a continuous function on a compact set $A \subset \mathbb{R}^N$. Verify that, setting:

$$F(x) \doteq \min_{y \in A} \left\{ F(y) + \frac{|x-y|}{\text{dist}(x, A)} - 1 \right\} \quad \text{for all } x \in \mathbb{R}^N \setminus A,$$

defines a continuous extension of $F$ on $\mathbb{R}^N$. This construction is due to Hausdorff and it provides a proof of the Tietze extension theorem.

(ii) Let $u : \mathbb{R}^N \to \mathbb{R}$ be a bounded, Borel function and let $\epsilon : \mathbb{R}^N \to (0, \infty)$ be continuous. Show that the function: $x \mapsto \mathcal{A}_{\epsilon(x)} u(x)$ is continuous on $\mathbb{R}^N$.

**Exercise 2.8**  Modify the construction of the ball walk to the *sphere walk* using the outline below.

(i) Let $\Omega_1 = \partial B_1(0) \subset \mathbb{R}^N$ and let $\mathbb{P}_1 = \sigma^{N-1}$ be the normalised spherical measure on the Borel $\sigma$-algebra $\mathcal{F}_1$ of subsets of $\Omega_1$ (see Example A.9). Define the induced probability spaces $(\Omega, \mathcal{F}, \mathbb{P})$ and $\{(\Omega_n, \mathcal{F}_n, \mathbb{P}_n)\}_{n=0}^\infty$ as in the case of the ball walk. For every $\epsilon \in (0, 1)$ and $x_0 \in \mathcal{D}$, let $\{X_n^{\epsilon, x_0} : \mathcal{D} \to \mathbb{R}^N\}_{n=0}^\infty$ be the sequence of random variables in:

$$X_0^{\epsilon, x_0} \equiv x_0 \quad \text{and for all } n \geq 1 \text{ and all } (w_1, \ldots, w_{n-1}) \in \Omega_{n-1}:$$

$$X_n^{\epsilon, x_0}(w_1, \ldots, w_n) = x_{n-1} + \left( \epsilon \wedge \frac{1}{2} \text{dist}(x_{n-1}, \partial \mathcal{D}) \right) w_n$$
$$\text{where} \quad x_{n-1} = X_{n-1}^{\epsilon, x_0}(w_1, \ldots, w_{n-1}).$$

Prove that $\{X_n^{\epsilon, x_0}\}_{n=0}^\infty$ is a martingale relative to the filtration $\{\mathcal{F}_n\}_{n=0}^\infty$ and that (2.3) holds for some random variable $X^{\epsilon, x_0} : \Omega \to \partial \mathcal{D}$.

(ii) For a continuous function $F : \partial \mathcal{D} \to \mathbb{R}$, define $u^\epsilon : \mathcal{D} \to \mathbb{R}$ according to (2.5). Show that $u^\epsilon$ is Borel-regular and that it satisfies:

$$u^\epsilon(x) = \fint_{\partial B_{\epsilon \wedge \frac{1}{2}\text{dist}(x, \partial \mathcal{D})}(x)} u^\epsilon(y) \, d\sigma^{N-1}(y) \quad \text{for all } x \in \mathcal{D}.$$

(iii) Deduce that if $F$ has a harmonic extension $u$ on $\bar{\mathcal{D}}$ as in (2.8), then $u^\epsilon = u$ for all $\epsilon \in (0, 1)$.

## 2.3 The ball walk and harmonic functions

The main result of this Section states that the uniform limits of values $\{u^\epsilon\}_{\epsilon \to 0}$ of the ball walk that we introduced in Section 2.2, are automatically harmonic. The proof relies on checking that each limiting function $u$ satisfies the mean value property on spheres. This is achieved by applying Doob's theorem to $u^\epsilon$ evaluated along its own walk process $\{X_n^{\epsilon, x_0}\}_{n=0}^\infty$, and choosing to stop on exiting the ball whose boundary coincides with the given sphere.

**Theorem 2.9**  *Let $J \subset (0, 1)$ be a sequence decreasing to $0$. Assume that $\{u^\epsilon\}_{\epsilon \in J}$ defined in (2.5), converges locally uniformly in $\mathcal{D}$, as $\epsilon \to 0$, $\epsilon \in J$, to some $u \in C(\mathcal{D})$. Then $u$ must be harmonic.*

*Proof*  **1.** In virtue of Theorem C.19, it suffices to prove that:

$$u(x_0) = \fint_{\partial B_r(x_0)} u(y) \, d\sigma^N(y) \qquad \text{for all } B_{2r}(x_0) \subset \mathcal{D}. \tag{2.9}$$

Fix $x_0 \in \mathcal{D}$ and $r \leq \frac{1}{2}\text{dist}(x_0, \partial\mathcal{D})$, and for each $\epsilon \in J$ consider the following random variable $\tau^\epsilon : \Omega \to \mathbb{N} \cup \{+\infty\}$:

$$\tau^\epsilon = \inf\{n \geq 1; \ X_n^{\epsilon, x_0} \notin B_r(x_0)\},$$

where $\{X_n^{\epsilon, x_0}\}_{n=0}^\infty$ is the usual sequence of the token positions (2.2) in the $\epsilon$-ball walk started at $x_0$. Clearly, $\tau^\epsilon$ is finite a.s. in view of convergence to the boundary in (2.3) and it is a stopping time relative to the filtration $\{\mathcal{F}_n\}_{n=0}^\infty$. By Corollary 2.5, Doob's theorem (Theorem A.31 (ii)) yields:

$$u^\epsilon(x_0) = \mathbb{E}[u^\epsilon \circ X_0^{\epsilon, x_0}] = \mathbb{E}[u^\epsilon \circ X_{\tau}^{\epsilon, x_0}],$$

while by passing to the limit with $\epsilon \to 0$ we obtain, by uniform convergence:

$$u(x_0) = \lim_{\epsilon \to 0, \epsilon \in J} \mathbb{E}[u^\epsilon \circ X_\tau^{\epsilon, x_0}] = \lim_{\epsilon \to 0, \epsilon \in J} \int_{B_{r+\epsilon}(x_0) \setminus B_r(x_0)} u(y) \, d\sigma_\epsilon(y). \tag{2.10}$$

The Borel probability measures $\{\sigma_\epsilon\}_{\epsilon \in (0, r)}$ are here defined on $\bar{B}_{2r}(x_0) \setminus B_r(x_0)$ by the push-forward procedure, as in Exercise A.8:

$$\sigma_\epsilon(A) \doteq \mathbb{P}(X_\tau^{\epsilon, x_0} \in A).$$

**2.** We now identify the limit in the right hand side of (2.10). Observe that, by construction, the measures $\sigma_\epsilon$ are rotationally invariant. Further, by Prohorov's theorem (Theorem A.10), each subsequence of $\{\sigma_\epsilon\}_{\epsilon \to 0, \epsilon \in J}$ has a further subsequence that converges (weakly-$*$) to a Borel probability measure $\mu$ on $\bar{B}_{2r}(x_0) \setminus B_r(x_0)$. Since each $\sigma_\epsilon$ is supported in $B_{r+\epsilon}(x_0) \setminus B_r(x_0)$, the limit $\mu$ must be supported on $\partial B_r(x_0)$. Also, $\mu$ is rotationally invariant in view of the same property of each $\sigma_\epsilon$. Consequently, $\mu = \sigma^{N-1}$ must be the uniquely defined, normalised spherical measure on $\partial B_r(x_0)$ (see Exercises 2.10 and 2.11). As the limit does not depend on the chosen subsequence of $J$, we conclude:

$$\lim_{\epsilon \to 0,\, \epsilon \in J} \int_{\bar{B}_{2r}(x_0) \setminus B_r(x_0)} u(y)\, \mathrm{d}\sigma_\epsilon(y) = \fint_{\partial B_r(x_0)} u(y)\, \mathrm{d}\sigma^{N-1}(y).$$

Together with (2.10), this establishes (2.9) as claimed. $\square$

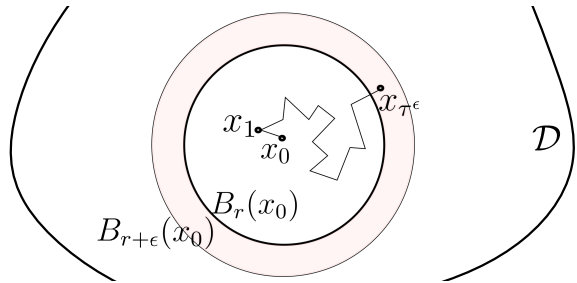

Figure 2.2 The stopping position $x_{\tau^\epsilon}$ in the proof of Theorem 2.9.

**Exercise 2.10** Show that every (weak-$*$) limit point of the family of probability measures $\{\sigma_\epsilon\}_{\epsilon \in (0,r)}$ defined in (2.10) must be rotationally invariant and supported on $\partial B_r(x_0)$.

**Exercise 2.11** Using the following outline, prove that the only Borel probability measure $\mu$ on $\partial B_1(0) \subset \mathbb{R}^N$ that is rotationally invariant, is the normalized spherical measure $\sigma^{N-1}$.

(i) Fix an open set $U \subset \partial B_1(0)$ and consider the sequence of Borel functions $\left\{ x \mapsto \frac{\mu(U \cap B(x, \frac{1}{n}))}{\mu(B(x, \frac{1}{n}))} \right\}_{n=1}^\infty$, where $B(x, r)$ denotes the $(N-1)$-dimensional curvilinear ball in $\partial B_1(0)$ centered at $x$ and with radius $r \in (0, 1)$. Apply Fatou's lemma (Theorem A.6) and Fubini's theorem (Theorem A.11) to the indicated sequence and deduce that:

$$\sigma^{N-1}(U) \leq \left( \liminf_{n \to \infty} \frac{\sigma^{N-1}(B(x, \frac{1}{n}))}{\mu(B(x, \frac{1}{n}))} \right) \cdot \mu(U), \qquad (2.11)$$

where both quantities $\sigma^{N-1}(B(x, \frac{1}{n}))$ and $\mu(B(x, \frac{1}{n}))$ are independent of $x \in \partial B_1(0)$ because of the rotational invariance.

(ii) Exchange the roles of $\mu$ and $\sigma^{N-1}$ in the above argument and conclude:

$$\lim_{n\to\infty} \frac{\sigma^{N-1}(B(x,\frac{1}{n}))}{\mu(B(x,\frac{1}{n}))} = 1.$$

Thus, $\mu(U) = \sigma^{N-1}(U)$ for all open sets $U$, so there must be $\mu = \sigma^{N-1}$.

This proof is due to Christensen (1970) and the statement above is a particular case of Haar's theorem on uniqueness of invariant measures on compact topological groups.

## 2.4 Convergence at the boundary and walk-regularity

We now investigate conditions assuring the validity of the uniform convergence assumption of Theorem 2.9. It turns out that such condition may be formulated independently of the boundary data $F$, only in terms of the behaviour of the ball walk (2.2) close to $\partial \mathcal{D}$, which is further guaranteed by a geometrical sufficient condition in the next Section. In Theorem 2.14 we will show how the boundary regularity of the process can be translated (via walk coupling) into the interior regularity, resulting in the existence of a harmonic extension $u$ of $F$ on $\bar{\mathcal{D}}$, and ultimately yielding $u^\epsilon = u$ for all $\epsilon \in (0,1)$, in virtue of Lemma 2.6.

**Definition 2.12** Consider the ball walk (2.2) on a domain $\mathcal{D} \subset \mathbb{R}^N$.

(a) We say that a boundary point $y_0 \in \partial \mathcal{D}$ is *walk-regular* if for every $\eta, \delta > 0$ there exists $\hat{\delta} \in (0, \delta)$ and $\hat{\epsilon} \in (0,1)$ such that:

$$\mathbb{P}(X^{\epsilon,x_0} \in B_\delta(y_0)) \geq 1 - \eta \quad \text{for all } \epsilon \in (0, \hat{\epsilon}) \text{ and all } x_0 \in B_{\hat{\delta}}(y_0) \cap \mathcal{D},$$

where $X^{\epsilon,x_0}$ is the limit in (2.3) of the $\epsilon$-ball walk started at $x_0$.
(b) We say that $\mathcal{D}$ is *walk-regular* if every $y_0 \in \partial \mathcal{D}$ is walk-regular.

**Lemma 2.13** *Assume that the boundary point $y_0 \in \partial \mathcal{D}$ of a given open, bounded, connected domain $\mathcal{D}$, is walk-regular. Then for every continuous $F : \partial \mathcal{D} \to \mathbb{R}$, the family $\{u^\epsilon\}_{\epsilon \to 0}$ defined in (2.5) satisfies the following. For every $\eta > 0$ there is $\hat{\delta} > 0$ and $\hat{\epsilon} \in (0,1)$ such that:*

$$|u^\epsilon(x_0) - F(y_0)| \leq \eta \quad \text{for all } \epsilon \in (0, \hat{\epsilon}) \text{ and all } x_0 \in B_{\hat{\delta}}(y_0) \cap \mathcal{D}. \quad (2.12)$$

*Proof* Given $\eta > 0$, let $\delta > 0$ satisfy:

$$|F(y) - F(y_0)| \leq \frac{\eta}{2} \quad \text{for all } y \in \partial \mathcal{D} \text{ such that } |y - y_0| < \delta.$$

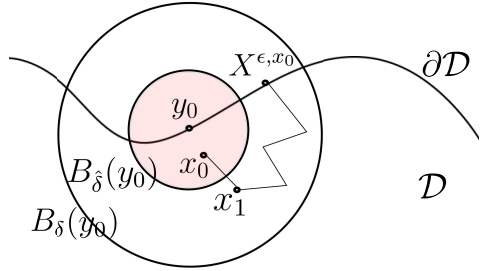

Figure 2.3 Walk-regularity of a boundary point $y_0 \in \partial \mathcal{D}$.

By Definition 2.12, we choose $\hat{\epsilon}$ and $\hat{\delta}$ corresponding to $\frac{\eta}{4\|F\|_\infty + 1}$ and $\delta$. Then:

$$|u^\epsilon(x_0) - F(y_0)| \leq \int_\Omega |F \circ X^{\epsilon,x_0} - F(y_0)| \, d\mathbb{P}$$

$$\leq \mathbb{P}(X^{\epsilon,x_0} \notin B_\delta(y_0)) \cdot 2\|F\|_\infty + \int_{X^{\epsilon,x_0} \in B_\delta(y_0)} |F \circ X^{\epsilon,x_0} - F(y_0)| \, d\mathbb{P}$$

$$\leq \frac{\eta}{4\|F\|_\infty + 1} \cdot 2\|F\|_\infty + \frac{\eta}{2} \leq \eta,$$

for all $x_0 \in B_{\hat{\delta}}(y_0) \cap \mathcal{D}$ and all $\epsilon \in (0, \hat{\epsilon})$. This completes the proof. □

By Lemma 2.13 and Theorem 2.9 we achieve the main result of this Chapter:

**Theorem 2.14** *Let $\mathcal{D}$ be walk-regular. Then, for every continuous $F : \partial \mathcal{D} \to \mathbb{R}$, the family $\{u^\epsilon\}_{\epsilon \in (0,1)}$ in (2.5) satisfies $u^\epsilon = u$, where $u \in C(\bar{\mathcal{D}})$ is the unique solution of the boundary value problem:*

$$\Delta u = 0 \quad \text{in } \mathcal{D}, \qquad u = F \quad \text{on } \partial \mathcal{D}.$$

*Proof* **1.** Let $F : \partial \mathcal{D} \to \mathbb{R}$ be a given continuous function. We will show that $\{u^\epsilon\}_{\epsilon \to 0}$ is "asymptotically equicontinuous in $\mathcal{D}$", i.e.: for every $\eta > 0$ there exists $\delta > 0$ and $\hat{\epsilon} \in (0, 1)$ such that:

$$|u^\epsilon(x_0) - u^\epsilon(y_0)| \leq \eta \qquad \text{for all } \epsilon \in (0, \hat{\epsilon}) \qquad (2.13)$$
$$\text{and all } x_0, y_0 \in \mathcal{D} \text{ with } |x_0 - y_0| \leq \delta.$$

Since $\{u^\epsilon\}_{\epsilon \to 0}$ is equibounded (by $\|F\|_\infty$), condition (2.13) imply that for every sequence $J \subset (0, 1)$ converging to $0$, one can extract a further subsequence of $\{u^\epsilon\}_{\epsilon \in J}$ that converges locally uniformly in $\bar{\mathcal{D}}$. Further, in view of (2.12) it follows that $u \in C(\bar{\mathcal{D}})$ and $u = F$ on $\partial \mathcal{D}$ (see Exercise 2.17). By Theorem 2.9, we get that $u$ is harmonic in $\mathcal{D}$ and the result follows in virtue of Lemma 2.6.

**2.** To show (2.13), fix $\eta > 0$ and choose $\bar\delta > 0$ such that $|F(y) - F(\bar y)| \leq \frac{\eta}{3}$ for all $y, \bar y \in \partial\mathcal{D}$ with $|y - \bar y| \leq 3\bar\delta$. By (2.12), for each $y_0 \in \partial\mathcal{D}$ there exists $\hat\delta(y_0) \in (0, \bar\delta)$ and $\hat\epsilon(y_0) \in (0, 1)$ satisfying:

$$|u^\epsilon(x_0) - F(y_0)| \leq \frac{\eta}{3} \qquad \text{for all } \epsilon \in (0, \hat\epsilon(y_0)) \text{ and all } x_0 \in B_{\hat\delta(y_0)}(y_0) \cap \mathcal{D}.$$

The family of balls $\{B_{\hat\delta(y)}(y)\}_{y \in \partial\mathcal{D}}$ is then a covering of the compact set $\partial\mathcal{D}$; let $\{B_{\hat\delta(y_i)}(y_i)\}_{i=1}^n$ be its finite sub-cover and set $\hat\epsilon = \min_{i=1\ldots n} \hat\epsilon(y_i)$. Clearly:

$$\partial\mathcal{D} + B_{2\delta}(0) \subset \bigcup_{i=1}^n B_{\hat\delta(y_i)}(y_i)$$

for some $\delta > 0$ where we additionally request that $\delta < \bar\delta$. This implies:

$$\begin{aligned}|u^\epsilon(x_0) - u^\epsilon(y_0)| \leq \eta \qquad &\text{for all } \epsilon \in (0, \hat\epsilon)\\ &\text{and all } x_0, y_0 \in (\partial\mathcal{D} + B_{2\delta}(0)) \cap \mathcal{D} \text{ with } |x_0 - y_0| \leq \delta.\end{aligned} \qquad (2.14)$$

**3.** To conclude the proof of (2.13), fix $\epsilon \in (0, \hat\epsilon \wedge \delta)$ and let $x_0, y_0 \in \mathcal{D}$ satisfy $\mathrm{dist}(x_0, \partial\mathcal{D}) \geq \delta$, $\mathrm{dist}(y_0, \partial\mathcal{D}) \geq \delta$ and $|x_0 - y_0| < \delta$. Define $\tau_\delta : \Omega \to \mathbb{N} \cup \{+\infty\}$:

$$\tau_\delta = \min\{n \geq 1;\ \mathrm{dist}(x_n, \partial\mathcal{D}) < \delta \text{ or } \mathrm{dist}(y_n, \partial\mathcal{D}) < \delta\},$$

where $\{x_n = X_n^{\epsilon, x_0}\}_{n=0}^\infty$ and $\{y_n = X_n^{\epsilon, y_0}\}_{n=0}^\infty$ denote the consecutive positions in the process (2.5) started at $x_0$ and $y_0$, respectively. It is clear that $\tau_\delta$ is finite $\mathbb{P}$-a.s. in view of convergence to the boundary in (2.3), and it is a stopping time relative to the filtration $\{\mathcal{F}_n\}_{n=0}^\infty$.

By Corollary 2.5 and Doob's theorem (Theorem A.31 (ii)) it follows that:

$$u^\epsilon(x_0) = \mathbb{E}[u^\epsilon \circ X_\tau^{\epsilon, x_0}] \quad \text{and} \quad u^\epsilon(y_0) = \mathbb{E}[u^\epsilon \circ X_\tau^{\epsilon, y_0}].$$

Since $|X_\tau^{\epsilon, x_0} - X_\tau^{\epsilon, y_0}| = |x_0 - y_0| < \delta$ and $X_\tau^{\epsilon, x_0}, X_\tau^{\epsilon, y_0} \in (\partial\mathcal{D} + B_{2\delta}(0)) \cap \mathcal{D}$ for a.e. $\omega \in \Omega$, we conclude by (2.14) that:

$$|u^\epsilon(x_0) - u^\epsilon(y_0)| \leq \int_\Omega |u^\epsilon \circ X_\tau^{\epsilon, x_0} - u^\epsilon \circ X_\tau^{\epsilon, y_0}|\, \mathrm{d}\mathbb{P} \leq \eta.$$

This ends the proof of (2.13) and of the Theorem. □

Walk-regularity is, in fact, equivalent to convergence of $u$ to the right boundary values. We have the following observation, converse to Lemma 2.13:

**Lemma 2.15** *If $y_0 \in \partial\mathcal{D}$ is not walk-regular, then there exists a continuous function $F : \partial\mathcal{D} \to \mathbb{R}$, such that for $u^\epsilon$ in (2.5) there holds:*

$$\limsup_{x \to y_0,\ \epsilon \to 0} u^\epsilon(x) \neq F(y_0).$$

*Proof* Define $F(y) = |y - y_0|$ for all $y \in \partial\mathcal{D}$. By assumption, there exists $\eta, \delta > 0$ and sequences $\{\epsilon_i\}_{i=1}^\infty$, $\{x_j \in \mathcal{D}\}_{j=1}^\infty$ such that:

$$\lim_{j\to\infty} \epsilon_j = 0, \quad \lim_{j\to\infty} x_j = y_0 \quad \text{and} \quad \mathbb{P}(X^{\epsilon_j,x_j} \notin B_\delta(y_0)) > \eta \quad \text{for all } j \geq 1,$$

where each $X^{\epsilon_j,x_j}$ above stands for the limiting random variable in (2.3) corresponding to the $\epsilon_j$- ball walk. By the nonnegativity of $F$, it follows that:

$$u^{\epsilon_j}(x_j) - F(y_0) = \int_\Omega F \circ X^{\epsilon_j,x_j} \, d\mathbb{P} \geq \int_{\{X^{\epsilon_j,x_j} \notin B_\delta(y_0)\}} F \circ X^{\epsilon_j,x_j} \, d\mathbb{P} > \eta\delta > 0,$$

proving the claim. $\square$

**Exercise 2.16** Show that if $\mathcal{D}$ is walk-regular then $\hat{\delta}$ and $\hat{\epsilon}$ in Definition 2.12 (a) can be chosen independently of $y_0$ (i.e. $\hat{\delta}$ and $\hat{\epsilon}$ depend only on the parameters $\eta$ and $\delta$).

**Exercise 2.17** Let $\{u_\epsilon\}_{\epsilon \in J}$ be an equibounded sequence of functions $u_\epsilon : \mathcal{D} \to \mathbb{R}$ defined on an open, bounded set $\mathcal{D} \subset \mathbb{R}^N$, and satisfying (2.12), (2.13) with some continuous $F : \partial\mathcal{D} \to \mathbb{R}$. Prove that $\{u_\epsilon\}_{\epsilon \in J}$ must have a subsequence that converges uniformly, as $\epsilon \to 0$, $\epsilon \in J$, to a continuous function $u : \bar{\mathcal{D}} \to \mathbb{R}$.

## 2.5 A sufficient condition for walk-regularity

In this Section we state a geometric condition (exterior cone condition) implying the validity of the walk-regularity condition introduced in Definition 2.12. We remark that the exterior cone condition in Theorem 2.19 may be weakened to the so-called exterior corkscrew condition, and that the analysis below is valid not only in the presently studied linear case of $p = 2$, but in the nonlinear setting of an arbitrary exponent $p \in (1, \infty)$ as well. This will be explained in Chapter 6, with proofs conceptually based on what follows.

We begin by observing a useful technical reformulation of the regularity condition in Definition 2.4. Namely, at walk-regular boundary points $y_0$ not only the limiting position of the ball walk may be guaranteed to stay close to $y_0$ with high probability, but the same local property may be, in fact, requested for the whole walk trajectory, with uniformly positive probability.

**Lemma 2.18** *Let $\mathcal{D} \subset \mathbb{R}^N$ be an open, bounded, connected domain. For a given boundary point $y_0 \in \partial\mathcal{D}$, assume that there exists $\theta_0 < 1$ such that for every $\delta > 0$ there are $\hat{\delta} \in (0, \delta)$ and $\hat{\epsilon} \in (0, 1)$ with the following property. For all $\epsilon \in (0, \hat{\epsilon})$ and all $x_0 \in B_{\hat{\delta}}(y_0) \cap \mathcal{D}$ there holds:*

$$\mathbb{P}(\exists n \geq 0 \; X_n^{\epsilon,x_0} \notin B_\delta(y_0)) \leq \theta_0, \tag{2.15}$$

where $\{X_n^{\epsilon,x_0}\}_{n=0}^{\infty}$ is the $\epsilon$-ball walk defined in (2.2). Then $y_0$ is walk-regular.

*Proof*  **1.** Fix $\eta, \delta > 0$ and let $m \in \mathbb{N}$ be such that:

$$\theta_0^m \leq \eta.$$

Define the tuples $\{\epsilon_k\}_{k=0}^m$, $\{\hat{\delta}_k\}_{k=0}^{m-1}$ and $\{\delta_k\}_{k=1}^m$ inductively, in:

$$\delta_m = \delta, \quad \epsilon_m = 1$$
$$\hat{\delta}_{k-1} \in (0, \delta_k), \quad \epsilon_{k-1} \in (0, \epsilon_k) \quad \text{for all} \ \ k = 1, \ldots, m \ \ \text{so that:}$$
$$\mathbb{P}(\exists n \geq 0 \ \ X_n^{x_0} \notin B_{\delta_k}(y_0)) \leq \theta_0 \qquad (2.16)$$
$$\text{for all} \ \ x_0 \in B_{\hat{\delta}_{k-1}}(y_0) \cap \mathcal{D} \ \ \text{and all} \ \ \epsilon \in (0, \epsilon_{k-1}),$$
$$\delta_{k-1} \in (0, \hat{\delta}_{k-1}) \quad \text{for all} \ \ k = 2, \ldots, m.$$

We finally set:

$$\hat{\epsilon} \doteq \epsilon_0 \wedge \min_{k=1,\ldots,m-1} |\hat{\delta}_k - \delta_k| \quad \text{and} \quad \hat{\delta} \doteq \hat{\delta}_0.$$

Fix $x_0 \in B_{\hat{\delta}}(y_0) \cap \mathcal{D}$ and $\epsilon \in (0, \hat{\epsilon})$. We will show that:

$$\mathbb{P}(\exists n \geq 0 \ \ X_n^{\epsilon,x_0} \notin B_{\delta_k}(y_0))$$
$$\leq \theta_0 \cdot \mathbb{P}(\exists n \geq 0 \ \ X_n^{\epsilon,x_0} \notin B_{\delta_{k-1}}(y_0)) \quad \text{for all} \ \ k = 2, \ldots, m. \qquad (2.17)$$

Together with the inequality in (2.16) for $k = 1$, the above bounds will yield:

$$\mathbb{P}(X^{\epsilon,x_0} \notin B_{2\delta}(y_0)) \leq \mathbb{P}(\exists n \geq 0 \ \ X_n^{\epsilon,x_0} \notin B_\delta(y_0)) \leq \theta_0^m \leq \eta.$$

Since $\eta$ and $\delta$ were arbitrary, the validity of the condition in Definition 2.12 will thus be justified, proving the walk-regularity of $y_0$.

**2.** Towards showing (2.17), we denote:

$$\tilde{\Omega} = \{\exists n \geq 0 \ \ X_n^{\epsilon,x_0} \notin B_{\delta_{k-1}}(y_0)\} \subset \Omega.$$

Without loss of generality, we may assume that $\mathbb{P}(\tilde{\Omega}) > 0$, because otherwise $\mathbb{P}(\exists n \geq 0 \ \ X_n^{\epsilon,x_0} \notin B_{\delta_k}(y_0)) \leq \mathbb{P}(\exists n \geq 0 \ \ X_n^{\epsilon,x_0} \notin B_{\delta_{k-1}}(y_0)) = 0$ and (2.17) holds then trivially. Consider the probability space $(\tilde{\Omega}, \tilde{\mathcal{F}}, \tilde{\mathbb{P}})$ defined by:

$$\tilde{\mathcal{F}} = \{A \cap \tilde{\Omega}; \ A \in \mathcal{F}\} \quad \text{and} \quad \tilde{\mathbb{P}}(A) = \frac{\mathbb{P}(A)}{\mathbb{P}(\tilde{\Omega})} \quad \text{for all} \ A \in \tilde{F}.$$

Also, let the measurable space $(\Omega_{fin}, \mathcal{F}_{fin})$ be given by: $\Omega_{fin} = \bigcup_{n=1}^{\infty} \Omega_n$ and by taking $\mathcal{F}_{fin}$ to be the smallest $\sigma$-algebra containing $\bigcup_{n=1}^{\infty} \mathcal{F}_n$. Then the following random variable $\tau_k : \tilde{\Omega} \to \mathbb{N}$:

$$\tau_k \doteq \min \{n \geq 1; \ X_n^{\epsilon,x_0} \notin B_{\delta_{k-1}}(y_0)$$

is a stopping time on $\tilde{\Omega}$ with respect to the induced filtration $\{\tilde{\mathcal{F}}_n = \{A \cap \tilde{\Omega}; \ A \in \mathcal{F}_n\}\}_{n=0}^\infty$. We consider two further random variables below:

$$Y_1 : \tilde{\Omega} \to \Omega_{fin} \qquad Y_1(\{w_i\}_{i=1}^\infty) \doteq \{w_i\}_{i=1}^{\tau_k}$$
$$Y_2 : \tilde{\Omega} \to \Omega \qquad Y_2(\{w_i\}_{i=1}^\infty) \doteq \{w_i\}_{i=\tau_k+1}^\infty$$

and observe that they are independent, namely:

$$\tilde{\mathbb{P}}(Y_1 \in A_1) \cdot \tilde{\mathbb{P}}(Y_2 \in A_2) = \tilde{\mathbb{P}}(\{Y_1 \in A_1\} \cap \{Y_2 \in A_2\})$$
$$\text{for all } A_1 \in \mathcal{F}_{fin}, \ A_2 \in \mathcal{F}.$$

We now apply Lemma A.21 to $Y_1, Y_2$ and to the indicator function:

$$Z(\{w_i\}_{i=1}^s, \{w_i\}_{i=s+1}^\infty) \doteq \mathbb{1}_{\left\{\exists n \geq 0 \ X_n^{\epsilon, x_0}(\{w_i\}_{i=1}^\infty) \notin B_{\delta_k}(y_0)\right\}}$$

that is a random variable on the measurable space $\Omega_{fin} \times \Omega$, equipped with the product $\sigma$-algebra of $\mathcal{F}_{fin}$ and $\mathcal{F}$. It follows that:

$$\mathbb{P}(\exists n \geq 0 \ X_n^{\epsilon, x_0} \notin B_{\delta_k}(y_0)) = \int_{\tilde{\Omega}} Z \circ (Y_1, Y_2) \, \mathrm{d}\tilde{\mathbb{P}} = \int_{\tilde{\Omega}} f(\omega_1) \, \mathrm{d}\tilde{\mathbb{P}}(\omega_1),$$

where for each $\omega_1 = \{w_i\}_{i=1}^\infty \in \tilde{\Omega}$ we have:

$$f(\omega_1) = \mathbb{P}\Big(\{\bar{w}_i\}_{i=1}^\infty \in \Omega; \ \exists n \geq 0 \ X_n^{\epsilon, x_0}(\{w_i\}_{i=1}^{\tau_k}, \{\bar{w}_i\}_{i=\tau_k+1}^\infty) \notin B_{\delta_k}(y_0)\Big)$$
$$= \mathbb{P}(\tilde{\Omega}) \cdot \mathbb{P}\Big(\exists n \geq 0 \ X_n^{\epsilon, x_{\tau_k}} \notin B_{\delta_k}(y_0)\Big) \leq \mathbb{P}(\tilde{\Omega}) \cdot \theta_0,$$

in view of $x_{\tau_k} \in B_{\hat{\delta}_k}(y_0)$ and the construction assumption (2.16). This ends the proof of (2.17) and of the lemma. □

The main result of this Section is a geometric sufficient condition for walk-regularity. When combined with Theorem 2.14, it implies that every continuous boundary data $F$ admits the unique harmonic extension to any Lipschitz domain $\mathcal{D}$. This extension automatically coincides with all process values $u^\epsilon$, regardless of the choice of the upper bound sampling radius $\epsilon \in (0, 1)$

**Theorem 2.19** *Let $\mathcal{D} \subset \mathbb{R}^N$ be open, bounded, connected and assume that $y_0 \in \partial \mathcal{D}$ satisfies the* exterior cone condition, *i.e. there exists a finite cone $C \subset \mathbb{R}^N \setminus \mathcal{D}$ with the tip at $y_0$. Then $y_0$ is walk-regular.*

*Proof* The exterior cone condition assures the existence of a constant $R > 0$ such that for all sufficiently small $\rho > 0$ there exists $z_0 \in C$ satisfying:

$$|z_0 - y_0| = \rho(1 + R) \quad \text{and} \quad B_{R\rho}(z_0) \subset C \subset \mathbb{R}^N \setminus \mathcal{D}. \tag{2.18}$$

Let $\delta > 0$ be, without loss of generality, sufficiently small and define $z_0 \in \mathbb{R}^N$ as in (2.18) with $\rho = \hat{\delta}$, where we set $\hat{\delta} = \frac{\delta}{4+2R}$. We will show that condition (2.15) holds for all $\epsilon \in (0, 1)$.

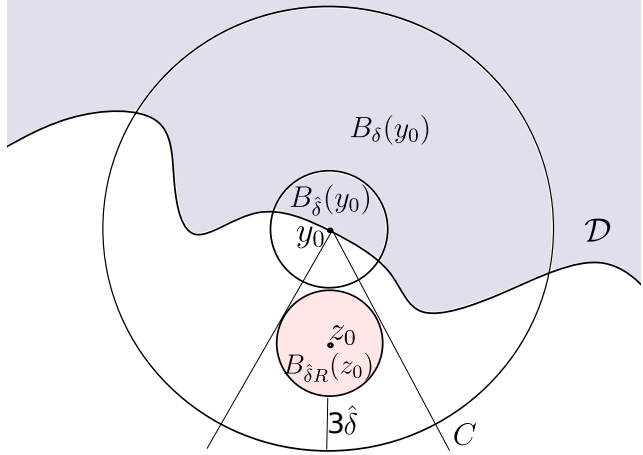

Figure 2.4 The concentric balls in the proof of Theorem 2.19.

Fix $x_0 \in B_{\hat{\delta}}(y_0) \cap \mathcal{D}$ and consider the profile function $v : (0, \infty) \to \mathbb{R}$ in:

$$v(t) = \begin{cases} \text{sgn}(N-2)\, t^{2-N} & \text{for } N \neq 2, \\ -\log t & \text{for } N = 2. \end{cases}$$

By Exercise C.27, the radial function $x \mapsto v(|x - z_0|)$ is harmonic in $\mathbb{R}^N \setminus \{z_0\}$, so in view of Lemma 2.6 the sequence of random variables $\{v \circ |X_n^{\epsilon,x_0} - z_0|\}_{n=0}^\infty$ is a martingale with respect to the filtration $\{\mathcal{F}_n\}_{n=0}^\infty$. Further, define the random variable $\tau : \Omega \to \mathbb{N} \cup \{+\infty\}$ by:

$$\tau \doteq \inf\{n \geq 0;\ X_n^{\epsilon,x_0} \notin B_\delta(y_0)\ ,$$

where we suppress the dependence on $\epsilon$ in the above notation. Applying Doob's theorem (Theorem A.31 (ii)) we obtain:

$$v(|x_0 - z_0|) = \mathbb{E}[v \circ |X_0^{\epsilon,x_0} - z_0|] = \mathbb{E}[v \circ |X_{\tau \wedge n}^{\epsilon,x_0} - z_0|] \qquad \text{for all }\ n \geq 0,$$

because for every $n \geq 0$, the a.s. finite random variable $\tau \wedge n$ is a stopping time. Passing to the limit with $n \to \infty$ and recalling the definition (2.3), now yields:

$$v(|x_0 - z_0|) = \int_{\{\tau < +\infty\}} v(|X_\tau^{\epsilon,x_0} - z_0|)\, d\mathbb{P} + \int_{\{\tau = +\infty\}} v(|X^{\epsilon,x_0} - z_0|)\, d\mathbb{P}.$$

Since $v$ is a decreasing function, this results in:

$$v((2+R)\hat{\delta}) \leq v(|x_0 - z_0|) \leq \mathbb{P}(\tau < +\infty) \cdot v((3+R)\hat{\delta}) + \mathbb{P}(\tau = +\infty) \cdot v(R\hat{\delta})$$
$$= \mathbb{P}(\tau < +\infty) \cdot \left(v((3+R)\hat{\delta}) - v(R\hat{\delta})\right) + v(R\hat{\delta}),$$

in view of the following bounds:

$$|x_0 - z_0| \leq |x_0 - y_0| + |y_0 - z_0| < (2+r)\hat{\delta},$$
$$|X_\tau^{\epsilon,x_0} - z_0| \geq |X_\tau^{\epsilon,x_0} - y_0| - |y_0 - z_0| = \delta - (1+R)\hat{\delta} = (3+R)\hat{\delta},$$
$$|X^{\epsilon,x_0} - z_0| \geq R\hat{\delta}.$$

Finally, noting that $v((3+R)\hat{\delta}) - v(R\hat{\delta}) < 0$ we obtain:

$$\mathbb{P}(\tau < +\infty) \leq \frac{v(R\hat{\delta}) - v((2+R)\hat{\delta})}{v(R\hat{\delta}) - v((3+R)\hat{\delta})} = \frac{v(R) - v(2+R)}{v(R) - v(3+R)}.$$

This establishes (2.15) with the constant $\theta_0 = \frac{v(R)-v(2+R)}{v(R)-v(3+R)} < 1$, that depends only on the dimension $N$ and the cone $C$. By Lemma 2.18, the proof is done. $\square$

**Remark 2.20** An alternative sufficient condition for walk-regularity is the simple-connectedness of $\mathcal{D} \subset \mathbb{R}^2$. The proof follows through the identification of $\{X_n^{\epsilon,x_0}\}_{n=0}^\infty$ as the discrete realisation of the Brownian path in Section 2.7* and applying the same reasoning as in the proof of Theorem 6.21. Indeed, in Chapter 6 we will give sufficient conditions for the so-called game-regularity, in the context of the Dirichlet problem for $p$-Laplacian, $p \in (1, \infty)$, encompassing and extending the classical discussion in the present Chapter.

## 2.6* The ball walk values and Perron solutions

In this Section we prove that all functions in the family $\{u^\epsilon\}_{\epsilon \in (0,1)}$ defined in (2.5) for a continuous $F : \partial \mathcal{D} \to \mathbb{R}$, are always one and the same function, coinciding with the so-called Perron solution of the Dirichlet problem:

$$\Delta u = 0 \quad \text{in } \mathcal{D}, \qquad u = F \quad \text{on } \partial \mathcal{D}. \qquad (2.19)$$

This material may be skipped at first reading, as it requires familiarity with more advanced PDE notions of Perron's method and Wiener's resolutivity. The related presentation in the general nonlinear case of $\Delta_p$, $p \in (1, \infty)$, can be found in Section C.7 of Appendix C. Below we recall this classical approach in the linear setting $p = 2$; for proofs we refer to the textbook by Helms (2014).

**Definition 2.21**  (i) A function $v \in C(\mathcal{D})$ is called *superharmonic* in $\mathcal{D}$, provided that for every $\bar{B}_r(x) \subset \mathcal{D}$ and every $h \in C(\bar{B}_r(x))$ that is harmonic in $B_r(x)$ and satisfies $h \leq v$ on $\partial B_r(x)$, there holds: $h \leq v$ in $B_r(x)$.

(ii) A function $v \in C(\mathcal{D})$ is *subharmonic* in $\mathcal{D}$, when $(-v)$ is superharmonic.

(iii) Given a continuous boundary data function $F : \partial \mathcal{D} \to \mathbb{R}$, we define the *upper and lower Perron solutions* to (2.19):

$$\bar{h}_F = \inf\{v \in C(\bar{\mathcal{D}}) \text{ superharmonic, such that } F \leq v \text{ on } \partial \mathcal{D}\},$$
$$\underline{h}_F = \sup\{v \in C(\bar{\mathcal{D}}) \text{ subharmonic, such that } v \leq F \text{ on } \partial \mathcal{D}\}.$$

The usual maximum principle argument implies that if $v_1, v_2 \in C(\bar{\mathcal{D}})$ are, respectively, subharmonic and superharmonic, and if $v_1 \leq v_2$ on $\partial \mathcal{D}$, then $v_1 \leq v_2$ in $\mathcal{D}$. In this comparison result, the conclusion may be in fact strengthened to: $v_1 < v_2$ or $v_1 \equiv v_2$ in $\mathcal{D}$. It follows that $\bar{h}_F$ and $\underline{h}_F$ are well defined functions, and also: $\underline{h}_F \leq \bar{h}_F$. One may further prove, by means of the harmonic lifting, that $\bar{h}_F$ and $\underline{h}_F$ are harmonic in $\mathcal{D}$. The celebrated *Wiener resolutivity theorem* in Wiener (1925) states the uniqueness of this construction:

**Theorem 2.22**  *Let $\mathcal{D} \subset \mathbb{R}^N$ be open, bounded and connected. Every boundary data $F \in C(\partial \mathcal{D})$ is resolutive, i.e. the two functions $\bar{h}_F$ and $\underline{h}_F$ coincide in $\mathcal{D}$. The resulting harmonic function is called the* Perron solution *to (2.19):*

$$h_F = \bar{h}_F = \underline{h}_F. \tag{2.20}$$

We remark that $h_F$ does not have to attain the prescribed boundary value $F(x)$ at each $x \in \partial \mathcal{D}$; it necessarily does so, however, for all points outside of a set whose 2-capacity is zero (see Section C.7).

By identifying the super- / subharmonic functions via mean value inequalities and comparing $u^\epsilon$ with $\bar{h}_F$ and $\underline{h}_F$, we obtain the main result of this Section:

---

**Theorem 2.23**  *Let $\mathcal{D} \subset \mathbb{R}^N$ be open, bounded, connected and let $F \in C(\partial \mathcal{D})$. For each $\epsilon \in (0, 1)$, functions $u^\epsilon$ in (2.5) satisfy: $u^\epsilon = h_F$ in $\mathcal{D}$.*

---

*Proof*  Let $v \in C(\bar{\mathcal{D}})$ be superharmonic and satisfy $F \leq v$ on $\partial \mathcal{D}$. Observe first that for any ball $\bar{B}_r(x) \subset \mathcal{D}$ we may apply Definition 2.21 to compare $v$ and the harmonic extension $u$ of $v_{|\partial B_r(x)}$ on $B_r(x)$ (see Exercise C.21) and get:

$$\fint_{\partial B_r(x)} v(y) \, d\sigma^{N-1}(y) = \fint_{\partial B_r(x)} u(y) \, d\sigma^{N-1}(y) = u(x) \leq v(x).$$

Integrating in polar coordinates, as in the proof of Theorem C.19, we obtain:

$$\fint_{B_r(x)} v(y)\,\mathrm{d}y = \frac{1}{|B_r(x)|} \int_0^r \int_{\partial B_s(x)} v(y)\,\mathrm{d}\sigma^{N-1}(y)\,\mathrm{d}s$$

$$\leq \frac{1}{|B_r(x)|} \int_0^r \int_{\partial B_s(x)} |\partial B_s(x)|\,\mathrm{d}s \cdot v(x) = v(x).$$

Fix $\epsilon \in (0,1)$ and $x_0 \in \mathcal{D}$. The sequence of random variables $\{v \circ X_n^{\epsilon,x_0}\}_{n=0}^\infty$ along the ball walk $\{X_n^{\epsilon,x_0}\}_{n=0}^\infty$ defined in (2.2), is then a supermartingale with respect to the filtration $\{\mathcal{F}_n\}_{n=0}^\infty$, because:

$$\mathbb{E}(v \circ X_n^{\epsilon,x_0} \mid \mathcal{F}_{n-1}) = \fint_{\epsilon \wedge \mathrm{dist}(X_{n-1},\partial\mathcal{D})} v(y)\,\mathrm{d}y \leq v \circ X_{n-1}^{\epsilon,x_0} \quad a.s.$$

Consequently: $\mathbb{E}[v \circ X_n^{\epsilon,x_0}] \leq \mathbb{E}[v \circ X_0] = v(x_0)$. Passing to the limit with $n \to \infty$ and recalling the boundary comparison assumption, finally yields:

$$v(x_0) \geq \mathbb{E}[v \circ X^{\epsilon,x_0}] \geq \mathbb{E}[F \circ X^{\epsilon,x_0}] = u^\epsilon(x_0).$$

We conclude that $\bar{h}_F \geq u^\epsilon$ by taking the infimum over all $v$ as above. Since by a symmetric argument: $\underline{h}_F \leq u^\epsilon$, the result follows in virtue of (2.20). □

## 2.7* The ball walk and Brownian trajectories

In this Section we show that the ball walk, introduced in Section 2.2, can be seen as a discrete realisation of the Brownian motion. In particular, we will deduce the same result as in Section 2.6*, namely that all functions in the family $\{u^\epsilon\}_{\epsilon \in (0,1)}$ in (2.5) are always one and the same function. This material may be skipped at first reading; it is slightly more advanced and necessitates familiarity with the construction of Brownian motion in Appendix B.

We start with some elementary technical observations. Denote $(\Omega_\mathcal{B}, \mathcal{F}_\mathcal{B}, \mathbb{P}_\mathcal{B})$ the probability space on which the standard $N$-dimensional Brownian motion $\{\mathcal{B}_t^N\}_{t\geq 0}$ is defined. We consider the product probability space $(\bar{\Omega}, \bar{\mathcal{F}}, \bar{\mathbb{P}}) = (\Omega_\mathcal{B}, \mathcal{F}_\mathcal{B}, \mathbb{P}_\mathcal{B}) \times (\Omega, \mathcal{F}, \mathbb{P})$ with the space $(\Omega, \mathcal{F}, \mathbb{P})$ in Section 2.2, and denote its elements by $(\omega_\mathcal{B}, \omega)$ with $\omega = \{w_i\}_{i=1}^\infty \in B_1(0)^\mathbb{N}$. Clearly, $\{\mathcal{B}_t^N\}_{t\geq 0}$ is also a standard Brownian motion on $(\bar{\Omega}, \bar{\mathcal{F}}, \bar{\mathbb{P}})$. We further denote the product $\sigma$-algebras $\bar{\mathcal{F}}_t = \mathcal{F}_t \times \mathcal{F}$, so that $\bar{\mathcal{F}}_s \subset \bar{\mathcal{F}}_t \subset \bar{\mathcal{F}}$ for all $0 \leq s \leq t$; for every $s \in [0,t]$ the random variable $\mathcal{B}_s^N$ is $\bar{\mathcal{F}}_t$-measurable.

We call $\bar{\tau} : \bar{\Omega} \to [0,\infty]$ a stopping time on $(\bar{\Omega}, \bar{\mathcal{F}}, \bar{\mathbb{P}})$ provided that $\{\bar{\tau} \leq t\} \in \bar{\mathcal{F}}_t$ for all $t \geq 0$ and that $\bar{\mathbb{P}}(\bar{\tau} = +\infty) = 0$. Then, the random variable $\mathcal{B}_{\bar{\tau}}^N$ is $\bar{\mathcal{F}}_{\bar{\tau}}$-measurable, namely: $\{\mathcal{B}_{\bar{\tau}}^N \in A\} \cap \{\bar{\tau} \leq t\} \in \bar{\mathcal{F}}_t$ for all Borel $A \subset \mathbb{R}^N$ and all $t \geq 0$, which can be proved as in Lemma B.21.

Let now $\mathcal{D} \subset \mathbb{R}^N$ be open, bounded, connected and fix $x_0 \in \mathcal{D}$, $\epsilon \in (0, 1)$. We inductively define the sequence of random variables $\bar{\tau}_k : \bar{\Omega} \to [0, \infty]$ in:

$$\bar{\tau}_0 = 0,$$
$$\bar{\tau}_{k+1}(\omega_\mathcal{B}, \{w_i\}_{i=1}^\infty)$$
$$= \min \left\{ t \geq \bar{\tau}_k; \; \left| \mathcal{B}_t^N(\omega_\mathcal{B}) - \mathcal{B}_{\bar{\tau}_k(\omega_\mathcal{B}, \omega)}^N(\omega_\mathcal{B}) \right| \right. \quad (2.21)$$
$$\left. = (\epsilon \wedge \mathrm{dist}(x_0 + \mathcal{B}_{\bar{\tau}_k(\omega_\mathcal{B}, \omega)}^N(\omega_\mathcal{B}), \partial \mathcal{D})) |w_{k+1}| \right\},$$

and also:

$$\bar{\tau}(\omega_\mathcal{B}, \omega) = \min \{ t \geq 0; \; x_0 + \mathcal{B}_t^N(\omega_B) \in \partial \mathcal{D} . \quad (2.22)$$

**Lemma 2.24** *Each $\bar{\tau}_k$ in (2.21) and $\bar{\tau}$ in (2.22) is a stopping time on $(\bar{\Omega}, \bar{\mathcal{F}}, \bar{\mathbb{P}})$. Moreover, $\bar{\tau}_k$ converge to $\bar{\tau}$ as $k \to \infty$, a.s. in $\bar{\Omega}$.*

□

Given a continuous boundary function $F : \partial \mathcal{D} \to \mathbb{R}$, recall that:

$$u(x_0) = \int_{\bar{\Omega}} F \circ (x_0 + \mathcal{B}_{\bar{\tau}}^N) \, d\bar{\mathbb{P}} \quad (2.23)$$

defines a harmonic function $u : \mathcal{D} \to \mathbb{R}$, in virtue of Corollary B.29 that builds on the classical construction and discussion of Brownian motion presented in Appendix B. As in Remark 2.3, we view $F$ as a restriction of some $F \in C(\bar{\mathcal{D}})$. Then, by Lemma 2.24 we also have:

$$u(x_0) = \lim_{k \to \infty} \int_{\bar{\Omega}} F \circ (x_0 + \mathcal{B}_{\bar{\tau}_k}^N) \, d\bar{\mathbb{P}}.$$

On the other hand, we recall that in (2.6) we defined:

$$u^\epsilon(x_0) = \int_\Omega F \circ X^{\epsilon, x_0} \, d\mathbb{P} = \lim_{k \to \infty} \int_\Omega F \circ X_k^{\epsilon, x_0} \, d\mathbb{P}.$$

**Theorem 2.25** *For all $\epsilon \in (0, 1)$ and all $x_0 \in \mathcal{D}$ there holds: $u^\epsilon(x_0) = u(x_0)$. In fact, we have:*

$$\mathbb{P}_\mathcal{B}(x_0 + \mathcal{B}_{\bar{\tau}}^N \in A) = \mathbb{P}(X^{\epsilon, x_0} \in A) \qquad \textit{for all Borel } A \subset \mathbb{R}^N. \quad (2.24)$$

□

**Exercise 2.26** Modify the arguments in this Section to the setting of the sphere walk introduced in Exercise 2.8. Follow the outline below:

(i) Given $x_0 \in \mathcal{D}$ and $\epsilon \in (0, 1)$, show that the following are stopping times on $(\Omega_\mathcal{B}, \mathcal{F}_\mathcal{B}, \mathbb{P}_\mathcal{B})$:

$$\tau_0 = 0, \qquad \tau_{k+1} = \min \left\{ t \geq \tau_k; \; |\mathcal{B}_t^N - \mathcal{B}_{\tau_k}^N| = \epsilon \wedge \frac{1}{2} \mathrm{dist}(x_0 + \mathcal{B}_{\tau_k}^N, \partial \mathcal{D}) \right\},$$

that converge a.s. as $k \to \infty$, to the exit time:

$$\tau = \min \{ t \geq 0; \; \mathcal{B}_t^N \in \partial \mathcal{D} - x_0 \}.$$

(ii) Let $(\Omega, \mathcal{F}, \mathbb{P})$ and $\{X_n^{\epsilon,x_0}\}_{n=0}^{\infty}$ be as in Exercise 2.8 (i), and define $u^\epsilon : \mathcal{D} \to \mathbb{R}$ according to (2.5) and (2.6). Prove that the push-forward of $\mathbb{P}$ on $\mathcal{D}$ via $X_k^{\epsilon,x_0}$, coincides with the push-forward of $\mathbb{P}_\mathcal{B}$ via $x_0 + \mathcal{B}_{\tau_k}^N$, for every $k \geq 0$. Consequently, $u^\epsilon(x_0) = \int_\Omega F \circ (x_0 + \mathcal{B}_\tau^N)\, d\mathbb{P}$, which is the harmonic extension of a given $F \in C(\partial \mathcal{D})$, independent of $\epsilon \in (0, 1)$.

## 2.8 Bibliographical notes

All constructions, statements of results and proofs in this Chapter have their continuous random process counterparts through Brownian motion, see Mörters and Peres (2010). The ball walk can be seen as a modification of the sphere walk in Exercise 2.8, which in turn is one of the most commonly used methods for sampling from harmonic measure, proposed in Muller (1956).

The definition of the walk-regularity of a boundary point $y_0$, which in the context of Section 2.7* can be rephrased as:

$$\forall \eta, \delta > 0 \quad \exists \hat{\delta} \in (0, \delta) \quad \forall x_0 \in B_{\hat{\delta}}(y_0) \cap \mathcal{D} \qquad \mathbb{P}(x_0 + \mathcal{B}_\tau^N \in B_\delta(y_0)) \geq 1 - \eta,$$

is equivalent to the classical definition given in Doob (1984):

$$\mathbb{P}_\mathcal{B}\Big(\inf\{t > 0;\ y_0 + \mathcal{B}_t^N \in \mathbb{R}^N \setminus \mathcal{D}\} = 0\Big) = 1;$$

this equivalence will be shown in Section 3.7*. The above property is further equivalent to the classical potential theory 2-regularity of $y_0$ in Definition C.46. Its equivalence with the Wiener regularity criterion, stating that $\mathbb{R}^N \setminus \mathcal{D}$ is 2-thick at $y_0$ (compare Definition C.47 (ii)) can be proved directly, see Mörters and Peres (2010) for a modern exposition. In working out the proofs of this Chapter and the analysis in Section 2.7*, the author has largely benefited from the aforementioned book and from personal communications with Y. Peres.

Various averaging principles and related random walks in the Heisenberg group were discussed in Lewicka et al. (2019). In papers by Lewicka and Peres (2019a,b), Laplace's equation augmented by the Robin boundary conditions has been studied from the viewpoint of the related averaging principles in $C^{1,1}$-regular domains. There, the asymptotic Hölder regularity of the values of the $\epsilon$-walk has been proved, for any Hölder exponent $\alpha \in (0, 1)$ and up to the boundary of $\mathcal{D}$, together with the interior asymptotic Lipschitz equicontinuity.

The "ellipsoid walk" linked to the elliptic problem: $\text{Trace}(A(x)\nabla^2 u(x)) = 0$ has been analyzed in Arroyo and Parviainen (2019). For bounded, measurable coefficients matrix $A$ satisfying $\det A = 1$, and uniformly elliptic with the elliptic distortion ratio that is close to 1 in $\mathcal{D}$, this lead to proving the local asymptotic uniform Hölder continuity of the associated process values $u^\epsilon$.